\documentclass[11pt,oneside]{amsart}

   \usepackage[english]{babel}
   \usepackage{color}
\usepackage{graphicx,color}

 \usepackage{amssymb}

\usepackage{geometry}
\usepackage{enumerate}
\usepackage{xcolor}
\usepackage{color}
\usepackage{multicol}

\usepackage[colorlinks=true, pdfstartview=FitV, linkcolor=blue, citecolor=blue, urlcolor=blue]{hyperref}

\usepackage[normalem]{ulem}
\usepackage{xcolor}
			\usepackage{amsmath,amscd}
			\usepackage{hyperref}
            \usepackage{fontenc}
            \usepackage{textcomp}
\usepackage{inputenc}
      \newtheorem{theorem}{Theorem}[section]
       \newtheorem{theo}{Theorem}      
       \newtheorem{lemma}[theorem]{Lemma}
       
          \newtheorem{example}[theorem]{Example}

      \newtheorem{prop}[theorem]{Proposition}

       \numberwithin{equation}{section}

            \newtheorem{remark}[theorem]{Remark}

      \makeatletter
      \def\@setcopyright{}
      \def\serieslogo@{}
      \makeatother
   
\begin{document}

%



   \author{Adrien Boyer and Jean-Martin Paoli}
   \address{Universit\'e de Paris Cité - Paris - France }
   \email{adrien.boyer@imj-prg.fr}
	\address{Universit\'e de Corte - Corsica}
	\email{paoli\_j@univ-corse.fr }

  
   \title{Riesz operators and $L^{p}$-boundary representations for hyperbolic groups}


   \begin{abstract}
  We investigate $L^{p}$-boundary representations of hyperbolic groups. We prove that such representations are irreducible if and only if the corresponding Riesz operators are injective.   \end{abstract}

   \subjclass{Primary 43A15, 43A90; Secondary 22D12}

   \keywords{$L^{p}$-boundary representations, $L^p$-Radial RD,  irreducibility, Bader-Muchnik ergodic theorems.}

  

   \date{\today}


   \maketitle


\section{Introduction}
Consider a non-elementary hyperbolic group $\Gamma$ endowed with an invariant metric $d$,  satisfying some regularity assumptions, acting by measure class preserving transformations on its Gromov boundary $(\partial \Gamma,\nu_{d})$ equipped with  $\nu_{d}$ the so-called Patterson-Sullivan measure associated with $d$. This yields to a one parameter family of isometric representations on $L^{p}$-spaces  denoted by $(\pi_{t},L^{p}(\partial \Gamma,\nu_{d}))$, defined for almost every $\xi\in \partial \Gamma$ as \begin{align}\label{representation}
[\pi_{t}(\gamma)v](\xi)=\bigg(\frac{d\gamma_{*}\nu_{d}}{d\nu_{d}}(\xi)\bigg)^{\frac{1}{2}+t}v(\gamma^{-1}\xi).
\end{align}
  The above representation satisfies
   $\|\pi_{t}(\gamma) v\|_{p}=\|v\|_{p}$ for all $v\in L^{p}(\partial \Gamma,\nu_{d})$ and for all $\gamma \in \Gamma$, where  $p$ is such that $1/p=1/2+t$ with $-\frac{1}{2}<t<\frac{1}{2}$. We call these representations \emph{$L^{p}$-boundary representations} of hyperbolic groups.\\

 \emph{The boundary representation} of hyperbolic groups is nothing but $\pi_{0}$ and has been intensively studied this last decade and might be seen, from a dynamical point of view, as a generalization of ergodicity see \cite{BM}, \cite{BM2}, \cite{Ga}, \cite{Boy2}, \cite{BoyMa},
\cite{BoyP}, \cite{BG}, \cite{BLP}, \cite{Fink}, \cite{KS} \cite{KS2}, and in \cite{Ca}.
 In the papers \cite{Bou}, \cite{Boy} and \cite{BPi} the representations (\ref{representation}) have already been studied but rather as representations on Hilbert spaces. In this paper we focus on $L^{p}$-spaces. We characterize the irreducibility of $L^{p}$-boundary representations $(\pi_{t},L^{p}(\partial \Gamma,\nu_{d}))$ with $1/p=1/2+t$ (where $-1/2<t<1/2$) thanks to an intertwining operator associated with the metric $d$, denoted by  $\mathcal{I}_{t}$  satisfying $\mathcal{I}_{t}\pi_{t}=\pi_{-t}\mathcal{I}_{t}$. We prove that this a bounded operator $\mathcal{I}_{t}:L^{p}(\partial \Gamma,\nu)\rightarrow L^{q}(\partial \Gamma,\nu)$ with $1/q=1/2-t$, defined only for $0<t<1/2$ (see  Subsection \ref{intertwiner}). It already appears in the context of  hyperbolic groups and Hilbert spaces in \cite{BPi}, \cite{GAG} and in CAT(-1) spaces \cite{Bou} .\\
  Our main result is the following:
  \begin{theo} \label{mainT}
   For all $-1/2<t<1/2$ and $1/p=1/2+t$, the $L^{p}$-boundary representations $(\pi_{t},L^{p}(\partial \Gamma,\nu_{d}))$ corresponding to $d$ are irreducible if and only if the intertwiner $\mathcal{I}_{|t|}$ is injective. 
  \end{theo}
  
  \begin{remark}
  In particular, we prove that the representation $(\pi_{t},L^{p}/\ker \mathcal{I}_{t})$ is irreducible for $1/p=1/2+t$ with $0<t<1/2.$
  \end{remark}
  We also deduce the following result in the context of rank one semisimple Lie groups. We do not know if the following theorem is present in the literature. 
 \begin{theo}\label{latt}
   Let $G$ be a rank one semisimple Lie group of non compact type. 
Let $\Gamma$ be a lattice in $G$. The $L^{p}$-boundary  representations of $\Gamma$ corresponding to $\nu$ the unique $K$-invariant probability measure on its Poisson-Furstenberg boundary $G/P$ are irreducible for all $1<p<+\infty$, where $K$ is the maximal compact subgroup and $P$ the minimal parabolic subgroup of $G$. 
 \end{theo}

  Indeed, we derive the above theorem from a generalizations of an ergodic theorem \`a la Bader-Muchnik for $L^{p}$-boundary representations, see Theorem \ref{BML2} in the following.
  \subsection*{Notation}
Endow $\Gamma$ with the length function corresponding  to $d$, $|\cdot|: \Gamma \rightarrow \mathbb{R}^{+}$  defined as $|\gamma |=d(1,\gamma )$ where  $1$ 
is the identity element and  $\gamma \in \Gamma.$\\
Let $S^{\Gamma}_{n,R}:=\{ \gamma \in \Gamma| nR\leq |\gamma|<(n+1)R\}$ for $R>0$ and let $|S^{\Gamma}_{n,R}|$ be the cardinal of $S^{\Gamma}_{n,R}$.\\
As in \cite{Boy}, we recall the definition of a spherical function associated with $\pi_{t}$. This  is  the matrix coefficient: \begin{align}
\phi_{t}:\gamma \in \Gamma \mapsto \langle \pi_{t}(\gamma)\textbf{1}_{\partial \Gamma}, \textbf{1}_{\partial \Gamma}\rangle \in \mathbb{R}^{+},
\end{align}
 where $
\textbf{1}_{\partial \Gamma}$ stands for the characteristic function of $\partial \Gamma$. \\
It will be convenient to also introduce the \emph{continuous} function  (see \cite[Lemma 3.2]{BPi})
\begin{align}\label{lafonction}
\sigma_{t}:\xi\in \partial \Gamma \mapsto \mathcal{I}_{t}(\textbf{1}_{\partial \Gamma})(\xi)\in \mathbb{R}^{+}\;\;\;\; (t>0).\end{align}  This is also a strictly positive function.
 Let $M_{\sigma^{-1}_{t}}$ be the corresponding multiplication operator, that is a bounded operator  acting on any $L^{p}(\partial \Gamma,\nu_{d})$ for any $p>1$.\\
We denoted by $\mathcal{R}_{t}=M_{\sigma^{-1}_{t}} \mathcal{I}_{t}$ the Riesz operator as in \cite{BPi}.

\subsection*{Convergence results}
We deduce the above theorems from  a theorem  {\it  \`a la Bader-Muchnik} for $L^{p}$-boundary representations of hyperbolic groups. Surprisingly, this kind of theorem known in the Hilbertian  context holds for $L^{p}$-spaces.
 \begin{theorem}\label{BML2}
For  $R>0$ large enough, there exists a sequence of  measures $\mu_{n,R}:\Gamma \rightarrow \mathbb{R}^{+}$,    supported on $S^{\Gamma}_{n,R}$, satisfying $\mu_{n,R}(\gamma)\leq C /|S^{\Gamma}_{n,R}|$ for some $C>0$ independent of $n$  such that
for all $0<t<1/2  $,  for all $f,g\in C(\Gamma \cup \partial \Gamma)$, for all $v\in L^{p}(\partial \Gamma,\nu_{o})$ and $w\in L^{q}(\partial \Gamma,\nu_{o})$:
$$\sum_{\gamma \in S^{\Gamma}_{n,R}}\mu_{n,R}(\gamma) f(\gamma   ) g(\gamma^{-1} ) \frac{\langle \pi_{t}(\gamma)v,w\rangle }{\phi_{t}(\gamma)}\to \langle  g_{|_{\partial \Gamma}}\mathcal{R}_{t}(v),\textbf{1}_{\partial \Gamma}\rangle \langle  f_{|_{\partial \Gamma}},w  \rangle, $$
as $n\to +\infty$.
\end{theorem}

 \begin{remark}
 If $\mu$ denotes a finitely supported random walk on a non-elementary hyperbolic group $\Gamma$ and $d:=d_{\mu}$ denotes the corresponding Green metric, then $(\Gamma,d_{\mu})$  satisfies the assumptions of  Theorem \ref{mainT} and \ref{BML2}.
 \end{remark}

\subsection*{Acknowledgement}
The first author would like to thank Ian Tice for useful discussions about Proposition \ref{weakschur} and Nigel Higson for comments on Theorem \ref{latt}.

  \subsection{Structure of the paper}
  Section \ref{sec2} contains preliminaries on $\delta$-hyperbolic spaces and hyperbolic groups, Patterson-Sullivan measures and equidistribution results, $L^{p}$-boundary representations as well as the definition of the intertwiner $\mathcal{I}_{t}$ for $t>0.$\\
  In Section \ref{interpolation}, we recall some basic facts in interpolation theory and some known results about spherical functions for hyperbolic groups. We prove that $\mathcal{I}_{t}$ is a bounded operator from $L^{p}(\partial \Gamma,\nu)$ to $L^{q}(\partial \Gamma,\nu)$ with $1/p=1/2+t$ and $1/q=1/2-t$ where $0<t<1/2$.\\
  Section \ref{section4} is devoted to proofs of Theorem \ref{BML2} and Theorem \ref{mainT}. In particular, a new tool we use is a $L^p$-version of Radial Property RD for $L^p$-boundary representations.\\
  Section \ref{section5} is a discussion about the case of  rank one globally symmetric spaces of non compact type and we provide a proof of Theorem \ref{latt}. 
  
  \section {Preliminaries on geometrical setting}\label{sec2}
\subsection{The geometrical setting and the regularity assumptions of the metric} A nice reference is \cite{BH}.\\
A metric space \((X,d)\) is said to be \emph{hyperbolic} if there exists $\delta\geq 0$  and a\footnote{if the condition holds for some \(o\) and \(\delta\), then it holds for any \(o\) and \(2\delta\)} basepoint \(o\in X\) such that for any \(x,y,z\in X\)  one has
\begin{equation}\label{hyp}
  (x,y)_{o}\geq \min\{ (x,z)_{o},(z,y)_{o}\}-\delta,
\end{equation}
where \((x,y)_{o}\) stands for the \emph{Gromov product} of \(x\) and \(y\) from \(o\), that is
\begin{equation}
  (x,y)_{o}=\frac{1}{2}(d(x,o)+d(y,o)-d(x,y)).
\end{equation}

We consider \emph{proper} hyperbolic metric space  (the closed balls are compact). 

A sequence $(a_{n})_{n\in  \mathbb{N}}$ in $X$ converges at infinity if $(a_{i},a_{j})_{o}\rightarrow +\infty$ as $i,j$ goes to $+\infty$. Set  $(a_n)\sim(b_n)\Leftrightarrow (a_i,b_j)_o\to \infty$ as $i,j\to \infty$ : it defines an equivalence relation and the set of equivalence classes  (that does not depent on the base point) is denoted by $\partial X$ and is called the Gromov boundary of $X$. The topology on $X$ naturally extends to  $\overline{X}:=X\cup \partial X$ so that $\overline{X}$ and $\partial X$ are compact sets. The formula 
\begin{equation}\label{gromovextended}
(\xi,\eta)_{o}:= \sup \liminf_{i,j}(a_{i},b_{j})_{o} 
\end{equation}
(where the supremum is taken over all sequences $(a_n), (b_n)$ which represent $\xi$ and $\eta$ respectively)
allows to extend the Gromov product on $\overline{ X}\times \overline{ X}$ but in a \emph{non continuous way} in general. Moreover the boundary $\partial X$ carries a family of \emph{visual metrics}, depending on \(d\) and a real parameter \(\epsilon > 0\) denoted from now by $d_{o,\epsilon}$.  The metric space  $(\partial X,d_{o,\epsilon})$ is a compact subspace of the bordification $\overline{X}:=\partial X \cup X$ (also compact) and the open ball  centered at $\xi$ of radius $r$ with respect to $d_{o,\epsilon}$ will be denoted by $B(\xi,r)$.\\

It turns out that in general, the Gromov product does not extend continuously to the bordification, see for example \cite[Example 3.16]{BH}. Following the authors of \cite{NS}, we say that a hyperbolic space $X$ is $\epsilon$-good, where $\epsilon>0$, if the following
two properties hold for each base point $o\in X$:
\begin{itemize}
\item The Gromov product $(\cdot,\cdot)_{o}$ on $X$ extends continuously to the bordification $X\cup \partial X$.

\item  The map $d_{o,\epsilon}:(\xi,\eta)\in \partial X\times \partial X \mapsto  e^{-\epsilon(\xi,\eta)_{o}}$ is a metric on  $\partial X$.
\end{itemize}

The classical theory of $\delta$-hyperbolic spaces works under the assumption that the spaces are geodesic but to guarantee that the Gromov product extends continuously to the boundary, that is if two sequences  $a_{n},b_{m} \in X \to \xi,\eta \in \partial X$, then the Gromov product satisfies $(a_{n},b_{m})_{o}\to (\xi,\eta)_{o}$,  we shall work under the assumption of roughly geodesic spaces. In particular the conformal relation on the boundary holds: for all $x,y\in X$ and for all $\xi,\eta \in \partial X:$
\begin{equation}\label{conform}
d_{y,\epsilon}(\xi,\eta)=e^{\frac{\epsilon}{2}  \big(\beta_{\xi}(x,y)+\beta_{\eta}(x,y)\big)}d_{x,\epsilon}(\xi,\eta),
\end{equation}
where the Busemann function $\beta_{\cdot}(\cdot,\cdot)$ is defined as $$(\xi,x,y)\in\partial X \times X \times X \mapsto \beta_{\xi}(x,y):= \lim_{n\to +\infty}d(x,a_{n}) -d(y,a_{n}),$$ where $(a_{n})_{n\in \mathbb{N}}$ represents $\xi.$
Recall that  for all $\xi\in \partial X$ and $x,y\in X$  
\begin{equation}\label{buse}
\beta_{\xi}(x,y)=-d(x,y)+2(\xi,y)_{x}.
\end{equation}
A metric space $(X,d)$ is roughly geodesic  if there exists $C=C_X>0$ so that for all $x,y\in X$ there exists a rough geodesic joining $x$ and $y$,  that is map $r:[a,b]\subset \mathbb{R}\rightarrow X$ with $r(a)=x$ and $r(b)=y$ such that 

\begin{equation}\label{roughgeo} |t-s|-C_X \leq d(r(t),r(s))\leq  |t-s|+C_X
\end{equation}
for all $t,s\in [a,b]$.

We say that two rough geodesic rays 
$r,r':[0,+\infty)\rightarrow X$ are equivalent if \\ $\sup_{t}d(r(t),r'(t))<+\infty$. We  write $\partial_{r} X$ for the set of equivalence classes of rough geodesic rays. When $(X,d)$ is a proper roughly geodesic space,  $\partial X$ and $\partial_{r} X$ coincide.

\subsection{Hyperbolic groups}
For an introduction to theory of hyperbolic groups we refer to \cite{Gr} and \cite{G}.\\
Recall that a group $\Gamma$ acts properly discontinuously on a proper metric space if for every compact sets $K,L\subset X$, the set $|\{\gamma\in \Gamma\;;\; \gamma K\cap L\neq \emptyset\}|<\infty$.    A {\it group $\Gamma$ is said to be hyperbolic} if it acts by isometries on some proper hyperbolic metric space $(X,d)$ such that $X/\Gamma$ is compact. A hyperbolic group is necessarily finitely generated (by \v{S}varc-Milnor's lemma). For such $\Gamma$,  any  finite set of generators $\Sigma$  gives rise to a Cayley graph  $({\mathcal G}(\Gamma, \Sigma),d_{\Sigma})$  whose set of vertices are the elements of $\Gamma$, linked by length-one edges if and only if they differ by an element of $\Sigma$.  Every {\it geodesic} hyperbolic metric space $(X,d)$ on which $\Gamma$ acts by isometries properly discontinuously with compact quotient is quasi-isometric to a Cayley graph of a hyperbolic group. If $\Gamma$ is a hyperbolic group
 endowed with a left invariant metric quasi-isometric to a word metric, it turns out that the metric space $(\Gamma,d)$ is a proper roughly geodesic $\delta$-hyperbolic metric space, see for example \cite[Section 3.1]{Ga}.\\
The limit set of $\Gamma$ denoted by $\Lambda_{\Gamma}$ is the set of  accumulation points in $\partial X$ of an (actually any) orbit. Namely $\Lambda_{\Gamma}:=\overline{\Gamma . o}\cap \partial X$, with the closure in $\overline{X}$. We say that $\Gamma$ is non-elementary if $|\Lambda_{\Gamma}|>2$ (and in this case, $|\Lambda_{\Gamma}|=\infty$). If $\Gamma$ is non-elementary and  if the action is cocompact then  $\Lambda_{\Gamma}=\partial X$.

 Eventually, note that a combination of results due to Blach\`ere, Ha\" issinsky and Matthieu \cite{BHM} and of Nica and \v{S}pakula \cite{NS} provides
\begin{theorem}
A hyperbolic group acts by isometries, properly discontinuously and cocompactly on a proper roughly geodesic  $\epsilon$-good $\delta$-hyperbolic space.
\end{theorem}

\subsection{To sum up}\label{class}

We assume that the metric space we are considering satisfies the following conditions:
\begin{itemize}
\item The metric space $(X,d)$ is $\delta$-hyperbolic.

\item The metric space $(X,d)$ proper.
\item The metric space $(X,d)$ is roughly geodesic.
\item The metric space $(X,d)$ $\epsilon$-good with some $\epsilon>0,$
\end{itemize}
and we let a non-elementary group $\Gamma$ act on $(X,d)$ under the following conditions:

\begin{itemize}
\item The action of $\Gamma$ is by isometries.
\item The action is properly discontinuous.
\item The action is cocompact.
\end{itemize}

In other words, the group $\Gamma$ is a non-elementary hyperbolic group and thus $\Gamma$ is infinite, discrete, countable and non-amenable.

\subsection{The Patterson -Sullivan measure}
Fix such $(X,d)$, pick an origin $o\in X$ and set $B(o,R)=\{ x\in X|d(o,x)<R\}$. 
Consider a family of visual metrics $(d_{x,\epsilon})_{x\in X}$ associated with a parameter $\epsilon$. The compact metric space $(\partial X,d_{o,\epsilon})$  admits a Hausdorff measure of dimension  \begin{equation}\label{DQ} D:={Q\over \epsilon} \end{equation} 
where 
\begin{equation}\label{volumegrowth}
Q=Q_{\Gamma,d}:=\limsup_{R\to +\infty}\frac{1}{R}\log |\Gamma.o\cap B(o,R)|,
\end{equation}
 is the critical exponent of $\Gamma$ (w.r.t. its action on $(X,d)$).
This $D$-Hausdorff measure is nonzero, finite, unique up to a constant, and denoted by  \(\nu_{o}\) when we normalize it to be a probability. The fundamental property we use is the Ahlfors regularity: the support of $\nu_{o}$ is in $\partial X$ and we say that $\nu_{o}$ is Ahlfors regular of dimension \(D\), if we have the following estimate for the volumes of balls: there exists $C>0$ so that for all $\xi \in \Lambda_{\Gamma}$ for all \(r \leq Diam (\partial X)\)
\begin{equation}\label{Ahlfors}
 C^{-1} r^{D}\leq  \nu_{o}(B(r,\xi)) \leq C r^{D}.
\end{equation}
The \emph{class} of measures $\nu_o$ is invariant under the action of $\Gamma$ and independent of the choice of \(\epsilon\). We refer to \cite{Pa}, \cite{Su}, \cite{BMo} and \cite{Co} for Patterson-Sullivan measures theory.

\subsection{Shadows and control of Busemann functions}
\subsubsection*{Upper Gromov bounded by above}
This assumption appears in the work of Connell and Muchnik in \cite{CM} as well as in the work of Garncarek on boundary unitary representations \cite{Ga}. We say that a space $X$ is \emph{upper gromov bounded by above} with respect to $o$, if there exists a constant $M>0$ such that for all $x\in X$ we have
 \begin{equation} \sup_{\xi \in \partial X}(\xi,x)_{o}\geq d(o,x)-M.
 \end{equation}
 
  Morally, this definition allows us to choose a point in the boundary playing the role of the forward endpoint of a geodesic starting at $o$ passing through $x$ in the context of simply connected Riemannian manifold of negative curvature. \\
  We denote by $\hat{x}_{o}$ a point in the boundary satisfying 
 
 \begin{equation} \label{endpoint}
 (\hat{x}_{o},x)_{o}\geq d(o,x)-M.
 \end{equation}
  In particular, every roughly geodesic metric spaces are upper Gromov bounded by above 
   (see for example \cite[Lemma 4.1]{Ga}).  

\subsubsection{Definition of shadows}
Let $(X,d)$ be a roughly geodesic, $\epsilon$-good, $\delta$-hyperbolic space. 
Let $r>0$ and a base point $o \in X$.
Define a shadow for any $x\in X$ denoted by $O_{r}(o,x)$ as 
\begin{equation}
O_{r}(o,x):=\{ \xi\in \partial X | (\xi,x)_{o}\geq d(x,o)-r\}.
\end{equation}

\begin{lemma}\label{ombre} Let $r>M+\delta$.  Then
$$B(\hat{x}_{o},e^{-\epsilon(d(o,x)-r+\delta)})\subset O_{r}(o,x) \subset B(\hat{x}_{o},e^{-\epsilon(d(x,o)-r-\delta)}). $$
\end{lemma}
\begin{proof}

Assume $r>M+\delta$.
For the left inclusion we have
\begin{align*}
(\xi,x)_{o}&\geq \min \{ (\xi,\hat{x}_{o})_{o},(\hat{x}_{o},x)\}-\delta\\
&\geq\min \{d(o,x)-r+\delta,d(o,x)-M \}-\delta \\
&=d(o,x)-r \\
&\geq d(o,x)-r.
\end{align*}
For the other inclusion
\begin{align*}
(\xi,\hat{x}_{o})_{o}&\geq \min \{ (\xi,x)_{o},(\hat{x}_{o},x)_{o}\}-\delta\\
&\geq \min \{ d(x,o)-r,d(o,x)-M \}-\delta\\
&\geq d(o,x)-r -\delta.
\end{align*}
\end{proof}
The above lemma combined with Ahlfors regularity of $\nu_{o}$ provides
\begin{lemma}\label{shadow}
There exists $C>0$   such that for any $x \in X$, and for $r>M+\delta$ $$C^{-1}e^{-Qd(o,x)}\leq\nu_{o}(O_{r}(o,x))\leq C e^{-Qd(o,x)}.$$

\end{lemma}

Here is a lemma dealing with a covering and the multiplicity of a covering by shadows of the boundary. 
\begin{lemma}\label{multiplicity}  

\begin{enumerate} We have the two following properties:
\item  \label{item1}For $R>0$ large enough, there exists $r>0$ such that $$\cup_{\gamma \in S^{\Gamma}_{n,R}} O_{r}(o,\gamma o)=\partial X.$$
\item \label{item2}For all $R,r>0$ large enough, there exists an integer $m$ such that for all $\xi \in \partial X$ we have for all $n\in \mathbb{N}$, 
$\sum_{\gamma \in S^{\Gamma}_{n,R}} \textbf{1}_{O_{r}(o,\gamma o)}(\xi)\leq m.$
\end{enumerate}
\end{lemma} 
\begin{proof}
Let $\kappa$ be the diameter of a relatively compact fundamental domain of the action of $\Gamma$ on $X$ containing $o$. Set \begin{equation}\label{choice}
R>3(C_{X}+\kappa)
\end{equation} where $C_{X}$ is the constant coming from the assumption \ref{roughgeo}.\\
We prove (\ref{item1}). 
Let $\xi\in \partial X$ and consider $r_{o}$ a roughly geodesic representing $\xi$. Define $z_{\xi}:=r_{o}(nR+R/2)$. Hence, $nR\leq d(o,z_{\xi})< (n+1)R.$ Since the action is cocompact, there exists $\gamma \in \Gamma$ such that $d(\gamma o,z_{\xi})\leq \kappa.$ 
The choice of (\ref{choice}) ensures $ nR \leq d(\gamma o,o)<(n+1)R.$ Therefore, 

\begin{align*}
(\xi, \gamma o)_{o}&\geq \min\{(\xi,z_{\xi})_{o} ,(z_{\xi},\gamma o)_{o}\}-\delta\\
&\geq \min\{(n+1)R-3C_{X},nR-\kappa\}-\delta\\
&\geq nR -\kappa-\delta\\
&\geq |\gamma|-R -\kappa-\delta,
\end{align*}
and thus $\xi \in O_{r}(o,\gamma o)$ with $r=R +\kappa+\delta$.\\
We now prove (\ref{item2}). Take $r>0$ and $R$ satisfying (\ref{choice}). For any $\gamma \in S^{\Gamma}_{n,R}$
and for all $\xi \in O_{r}(o,\gamma o)$ 
\begin{align*}
(\gamma o,z_{\xi})_{o}&\geq \min\{(\gamma o,\xi)_{o},(\xi,z_{\xi})_{o} \}-\delta\\
&\geq \min\{|\gamma|-r ,(n+1)R-3C_{X} \}-\delta\\
&\geq \min\{nR-r ,(n+1)R-3C_{X} \}-\delta\\
&\geq \min\{nR-r ,nR-C_{X} \}-\delta\\
&\geq nR-r-\delta.
\end{align*}
By definition of the Gromov product we deduce that $d(\gamma o,z_{\xi})\leq \rho$ where $\rho=2(R+r+\delta)$.
Thus the set $\{\gamma \in S^{\Gamma}_{n,R} | O_{r}(o,\gamma o)\ni \xi \}$ is contained in $B(z_{\xi},\rho)\cap \Gamma .o$. Since the action is cocompact, by taking a positive constant $R'$ bigger than $\kappa>0$, we obtain $|B(z_{\xi},\rho)\cap \Gamma .o| \leq |B(o,R')\cap \Gamma. o| $ with $R'=\rho-\kappa$. Set $m:=|B(o,R')\cap \Gamma .o| $ to conclude the proof.
\end{proof}
Recall that there exists $M>0$, such that if $\gamma$ is an element of $\Gamma$, one can choose $\hat{\gamma}_{o}$ a point in $\partial X$ satisfying (\ref{endpoint})
 
 $$ (\hat{\gamma}_{o},\gamma o)_{o}\geq |\gamma|-M.$$

\begin{lemma}\label{crucial}
There exists $C>0$ such that for all $\xi \in \partial X$ there exists $g_{\xi}\in S^{\Gamma}_{n,R}$ such that for all $\gamma \in S^{\Gamma}_{n,R}$, 
$\beta_{\xi}(o,\gamma o)\leq  \beta_{\hat{\gamma}_{o}}(o,g_{\xi} o)+C.$
\end{lemma}
\begin{proof}
 Pick a point $\xi \in \partial X$.  Consider a roughly geodesic $r_{o}$ starting at $o$ representing $\xi$ and choose a point $z_{\xi}$ on it such that $nR\leq d(z_{\xi},o)< (n+1)R$ (since $R$ is large enough).
 We have for $\gamma \in S^{\Gamma}_{n,R}$  \begin{align*}
 (\gamma o, z_{\xi})_{o}&\geq \min\{(\gamma o, \hat{\gamma o})_{o}, (\hat{\gamma o}, z_{\xi})_{o}\}-\delta\\
 & \geq \min\{|\gamma|-M, (\hat{\gamma o}, z_{\xi})_{o}\}-\delta,
\end{align*}
and therefore either $(\hat{\gamma o}, z_{\xi})_{o}\leq  (\gamma o, z_{\xi})_{o} +\delta $ or $ (\gamma o, z_{\xi})_{o}\geq |\gamma|-M-\delta \geq d(o,z_{\xi})-M-R-\delta\geq (\hat{\gamma o}, z_{\xi})_{o}-M-R-\delta.$ In other words, 
\begin{equation}\label{equ}
(\hat{\gamma o}, z_{\xi})_{o}\leq  (\gamma o, z_{\xi})_{o} +C_{\delta,M,R}
\end{equation}
 with $C_{\delta,M,R}=M+R+\delta.$
It follows that
  \begin{align*}
 (\xi,\gamma o)_{o}&\geq \min\{ (\xi,z_{\xi})_{o}, (z_{\xi},\gamma o)_{o} \}-\delta \\
 &\geq \min\{ nR-C_{X} , (z_{\xi},\gamma o)_{o} \}-\delta\\
 &\geq \min\{ nR-C_{X} , (\hat{\gamma o}, z_{\xi})_{o}-C_{\delta,M,R} \}-\delta\\
 &\geq (\hat{\gamma o}, z_{\xi})_{o}-C',
\end{align*}
with $C'=C_{X}+R+\delta+C_{\delta,M,R}>0$.

Equality (\ref{buse}) implies that 
$\beta_{\xi}(o,\gamma o)\leq \beta_{\hat{\gamma o}}(o,z_{\xi})+2C'.$ Now, choose $R$ large enough so that there exists $g_{\xi}\in \Gamma$ with $d(g_{\xi}o,z_{\xi})\leq \kappa$ and $go \in S^{\Gamma}_{n,R}$, where $\kappa$ is the diameter of a fundamental domain of the action of $\Gamma$ on $X$ containing $o$. Then\\
$$\beta_{\hat{\gamma o}}(o,z_{\xi})=\beta_{\hat{\gamma o}}(o,g_{\xi} o)+\beta_{\hat{\gamma o}}(g_{\xi} o,z_{\xi})\leq \beta_{\hat{\gamma o}}(o,g_{\xi} o)+d(g_{\xi} o,z_{\xi})\leq \beta_{\hat{\gamma o}}(o,g_{\xi} o)+\kappa.$$
To conclude the proof write
$$\beta_{\xi}(o,\gamma o)\leq \beta_{\hat{\gamma o}}(o,g_{\xi}o)+2C'+\kappa.$$
\end{proof}

\subsection{Equidistribution \`a la Roblin-Margulis}\label{equid}The following theorem appears under this form for the first time in \cite[Theorem 3.2]{BG} and has been inspired by results in \cite{Ro} and \cite{Ma}. 
The unit Dirac mass centered at  $x\in X$ is denoted by $D_{x}$. 

\begin{theorem}\label{equi}
For any  $R>0$ large enough, there exists a sequence of  measures $\mu_{n,R}:\Gamma \rightarrow \mathbb{R}^{+}$ such that
\begin{enumerate}
 \item \label{growth}There exists $C>0$ satisfying  for all $n \in \mathbb{N}$ and all $\gamma \in S^{\Gamma}_{n,R}$ that $$\mu_{n,R}(\gamma)\leq C /  |S^{\Gamma}_{n,R}|.$$
 \item 
  We have the following convergence: $$\sum_{\gamma \in S^{\Gamma}_{n,R}}\mu_{n,R}(\gamma) D_{\gamma o} \otimes D_{\gamma^{-1} o} \rightharpoonup  \nu_{o}\otimes \nu_{o},$$ 
as $n\to +\infty$, for the weak* convergence in $C(\overline{X}\times \overline{X})$.
\end{enumerate}

\end{theorem}

\subsection{$L^{p}$-representations}
The expression  (\ref{representation}) of  $\pi_{t}$ defines an isometric $L^{p}$-representation of $\Gamma$  for the exponent
\begin{equation}
p=\frac{2}{1+2t},
\end{equation}
 with $0<t<1/2$.
Denote its conjugate exponent 
\begin{equation}
q=\frac{2}{1-2t}.
\end{equation}
Observe that the contragredient representation of $(\pi_{t},L^{p})$ is $(\pi_{-t},\overline{L^{q}})$
with respect to the (non-degenerate) pairing $$\langle\cdot, \cdot\rangle:(v,w)\in L^{p}\times \overline{L^{q}} \rightarrow   \int_{\partial X}v(\xi)\overline{w}(\xi) d\nu_{o}(\xi) \in \mathbb{C}.$$  

In particular, the adjoint operator  $\pi^{*}_{t}(\gamma)$  of $\pi_{t}(\gamma)$ is given for any $\gamma \in \Gamma$ by
\begin{equation}
\pi^{*}_{t}(\gamma)=\pi_{-t}(\gamma^{-1}).
\end{equation}
\subsection{An Intertwiner}\label{sec14}

Following \cite{BPi}, recall the definition  of the operator $\mathcal{I}_{t}$ for $t>0$: for almost every $\xi \in \partial X$ 
 \begin{equation}\label{intertwiner}
\mathcal{I}_{t}(v)(\xi):=\int_{\partial X} \frac{v(\eta)}{d^{(1-2t)D}_{o,\epsilon}(\xi,\eta)}d\nu_{o}(\eta).
\end{equation}
 It has already been observed in \cite{BPi} that $\mathcal{I}_{t}$ is a self-adjoint compact operator on $L^{2}(\partial X,\nu_{o})$. We will show that for $0<t<1/2$, the formula (\ref{intertwiner}) defines $\mathcal{I}_{t}$ as a bounded operator from $L^{p}(\partial X,\nu_{o})$ to $L^{q}(\partial X,\nu_{o})$ with with $1/p=1/2+t$ and $1/q=1/2-t$, see Proposition \ref{cont}. Moreover, working under the assumptions of $\epsilon$-good spaces guarantees that the operator  $\mathcal{I}_{t}$ intertwines $\pi_{t}$ and $\pi_{-t}$ from $L^{p}(\partial X,\nu_{o})$ to $L^{q}(\partial X,\nu_{o})$ thanks to the relation (\ref{conform}), see \cite[Proposition 3.17]{BPi}. Namely, for all $\gamma \in \Gamma$ and for all $v\in L^{p}(\partial X,\nu_{o}):$

\begin{equation}\label{intert}
\mathcal{I}_{t}\pi_{t}(\gamma)v=\pi_{-t}(\gamma)\mathcal{I}_{t}v.
\end{equation}

It will be also useful to consider \begin{equation}\label{sigmat}
\tilde{\sigma_{t}}:x\in \overline{X} \mapsto \int_{\partial X} e^{(1-2t)Q( x,\eta)_{o}}d\nu_{o}(\eta)\in \mathbb{R}^{+}.
\end{equation}
Observe that $\tilde{\sigma_{t}}$ restricted to $\partial X$ is nothing but $\sigma_{t}$ defined in (\ref{lafonction}).
We recall that
 the function $\tilde{\sigma_{t}}$ si continuous on $\overline{X},$ see \cite[Proposition 3.4 ]{BPi}.

\section{Interpolation theory: Strong inequality of type $(p,q)$ for the intertwining operator and Application of Riesz-Thorin Theorem}\label{interpolation}

The aims of this section is to provide some material of interpolation theory to prove the main result concerning the operator $\mathcal{I}_{t}$, with $0<t<1/2$ based on the \emph{weak-type Schur's test}. The connection of interpolation theory and Lorentz psaces with boundary representations are already in the paper \cite{Cow} and \cite[Chapter 6]{FiPi}. Note also the very recent work \cite{GAG}.

\subsection{Lorentz spaces, interpolation and applications}
We follow \cite{Cow}.
Let $(\Omega,\mu)$ be a measure space. If $f:\Omega \rightarrow \mathbb{C}$ is a measurable function then define the nonincreasing rearrangement of $f$ $$f^*(t):= \inf \bigg\{ s>0  |\mu\big( \{|f|>s \}\big)\leq t \bigg\}.$$
The real function $f^*$ is a positive nonincreasing function, equimeasurable with $|f|$ and right continuous. Define then the norm,  if $1<p,q<\infty$ as $$\|f\|_{L^{p,q}}=\bigg(\frac{p}{q} \int^{\infty}_{0}(t^{1/p}f^{*}(t))^{q}\frac{dt}{t}\bigg)^{1/q},$$ and $1<p<\infty$ with $q=\infty$:

$$\|f\|_{L^{p,\infty}}=\sup\{ t^{1/p}f^{*}(t)|t\in [0,+\infty]\}.$$
Define the Lorentz spaces for $1<p<+\infty$ and $1< q\leq +\infty$.
$$L^{p,q}(\Omega):=\{ f:\Omega \rightarrow \mathbb{C} \mbox{ measurable }| \|f\|_{p,q}<+\infty\}.$$
Here are some useful facts:
\begin{enumerate}
\item $L^{p,p}=L^{p}.$
\item $(L^{p,q})^{*}=L^{p',q'}.$
\item \label{point3} $L^{p,q_{1}}\subset L^{p,q_{2}}$ with $q_{1}<q_{2}$, that is $$\|v\|_{L^{p,q_{2}}}\leq C_{p,q_{1},q_{2}} \|v\|_{L^{p,q_{1}}},$$for some $C_{p,q_{1},q_{2}}>0.$
\end{enumerate}
Here is the \emph{} the fundamental tool, called \emph{the weak-type Schur's test} coming form interpolation theory.

\begin{prop}\label{weakschur}
Let $(X,\mathcal{M},\mu)$ and $(Y,\mathcal{N},\nu)$ be $\sigma$-finite measure spaces and let $1<p,q,r<\infty$ be such that $$\frac{1}{p}+\frac{1}{r}=\frac{1}{q}+1.$$
Let $k:X\times Y \mapsto \mathbb{R}$ be a measurable function and suppose that there exists $A>0$ such that
\begin{align*}
\|k(\cdot,y)\|_{L^{r,\infty}}&\leq A \mbox{ for a.e. } y \in Y, \\
\|k(x,\cdot)\|_{L^{r,\infty}}&\leq A \mbox{ for a.e.  } x \in X. 
\end{align*}
Therefore the formula $Tv(x)=\int_{Y}k(x,y)v(y)d\nu(y)$ defined a.e a function in $L^{q}(X,\mu)$ whenever $v$ is in $L^{p}(Y,\nu)$. And moreover for all $1\leq s \leq \infty $ there exists a constant $C=C_{p,q,r,s}>0$ depending on $p,q,r,s$ such that

$$\|Tv\|_{L^{q,s}}\leq C A\|v\|_{L^{p,s}}.$$
\end{prop} 

We refer to \cite[Proposition 6.1]{Tao} for a proof.

\subsubsection{Analogs of homogenous functions on the boundary}
\begin{lemma} Let $1<1/p<\infty.$ 
We have for all $\xi\in \partial X$ $$\frac{1}{d^{ pD}_{o,\epsilon} (\xi,\cdot)}\in L^{1/p,\infty}.$$ 
\end{lemma}
\begin{proof}
Let $s>0.$ We have for all $\xi \in \partial X$
 $$\nu_{o}(\{ \eta |d^{-pD}_{o,\epsilon}(\xi,\eta)>s\})=\nu_{o}(\{ \eta |d_{o,\epsilon}(\xi,\eta)<s^{-1/pD}\})=\nu_{o}\big(B(\xi,s^{-1/pD})\big).$$
We obtain for all $\xi \in \partial X$ and for $t>0:$
\begin{align*}
\bigg(\frac{1}{d^{ pD}_{o,\epsilon} (\xi,\cdot)}\bigg)^{*}(t)&=\inf \{ s>0| \nu_{o}(\{ \eta |d^{-pD}_{o,\epsilon}(\xi,\eta)>s)\leq t\}\\
&= \inf \{ s>0 | \nu_{o}\big(B(\xi,s^{-1/pD})\big) \leq t \}  \\
&\leq  \inf \{ s>0 |  C s^{-1/p} \leq t \}  \\
&= C^{p}t^{-p},
\end{align*}
for a constant $C>0$ coming from the Ahlfors regularity property \ref{Ahlfors}. It follows that $1/d^{ pD}_{o,\epsilon} (\xi,\cdot) \in L^{1/p,\infty}.$ 
\end{proof}

For $t>0$ we have for all $\xi\in \partial X$ $$\frac{1}{d^{(1-2t) D}_{o,\epsilon} (\xi,\cdot)}\in L^{r,\infty},$$
with $r=1/(1-2t).$ By symmetry, we have for all $\eta \in \partial X$ that $$\frac{1}{d^{(1-2t) D}_{o,\epsilon} (\cdot,\eta)}\in L^{r,\infty}$$ as well. We obtain the following result:
\begin{prop}\label{cont}
Let $0<t<1/2$ and let $p,q$ such that $1/p=1/2+t$ and $1/q=1/2-t$.
The operator $\mathcal{I}_{t}$ is bounded from $L^{p}(\partial X,\nu_{o})$ to $L^{q}(\partial X,\nu_{o})$.
\end{prop}
\begin{proof}
Note that $1/p+(1-2t)=1/2+t+(1-2t)=1+1/2-t=1+1/q.$ Proposition \ref{weakschur} implies that $$\|\mathcal{I}_{t}(v)\|_{L^{q,s}}\leq C_{t} \|v\|_{L^{p,s}} $$ for all $1\leq s\leq \infty$.
Pick $p \leq s=2/1-2t=q $. Thus $\|\mathcal{I}_{t}(v)\|_{L^{q,s}}=\|\mathcal{I}_{t}(v)\|_{L^{q}}$ and $\|v\|_{L^{p,s}}=\|v\|_{L^{p,q}}\leq C_{p,q}\|v\|_{L^{p,p}}=\|v\|_{L^{p}} $ by (\ref{point3}) to complete the proof.
\end{proof}

    \subsection{Consequences of spectral gap estimates and Riesz-Thorin theorem}\label{section6}
  
   \subsubsection{Spherical functions on hyperbolic groups }
 
As in \cite{Boy}, we recall the definition of a spherical function associated with $\pi_{t}$. This  is  the matrix coefficient: \begin{align}
\phi_{t}:\gamma \in \Gamma \mapsto \langle \pi_{t}(\gamma)\textbf{1}_{\partial X}, \textbf{1}_{\partial X}\rangle \in \mathbb{R}^{+},
\end{align}
and introduce the function $\omega_{t}(\cdot)$ for $t\in \mathbb{R}^{*} $ defined as 
\begin{equation}
\omega_{t}(x) =\frac{2 \sinh\big( tQ x\big) }{e^{2tQ}-1}.
\end{equation}
Note that $\omega_{t}(\cdot)$ is a positive function for all $t\in \mathbb{R}^{*} $ converging uniformly on compact sets of $\mathbb{R}$ to 
    $ \omega_{0}(x)=x $ as $t\to 0$, and $\omega_{-t}(x)=e^{2tQ}\omega_{t}(x)$.\\

In \cite{Boy},  the following estimates, called \emph{ Harish-Chandra Anker estimates}, naming related to  \cite{Ank} have been proved.
There exists $C>0$, such that for any $t\in \mathbb{R}$, we have for all $\gamma \in \Gamma$
    \begin{align}\label{HCHestimates}
C^{-1}e^{-\frac{1}{2}Q|\gamma| }\big(1+\omega_{|t|}(|\gamma|)\big)
 \leq   \phi_{t}(\gamma)\leq 
 Ce^{-\frac{1}{2}Q|\gamma|}\big(1+\omega_{|t|}(|\gamma|)\big).
   \end{align}

Set for all $x\in \mathbb{R}$ \begin{equation}
\widetilde{\phi_{t}}(x):= e^{-\frac{1}{2}Qx}\big(1+\omega_{|t|}(x)\big).
\end{equation}

\subsubsection{A $L^{p}$-spectral inequality} We briefly recall some facts.
In \cite{Boy} the following spectral inequality, generalizing the so called ``Haagerup Property" or Property RD has been proved.  Pick $R>0$ large enough. There exists $C>0$ such that  for any $t \in \mathbb{R}$ and  for all $f_{n}\in \mathbb{C}[\Gamma]$ supported in $S^{\Gamma}_{n,R}$, we have 
\begin{equation}\label{RDgeneral}
\|\pi_{t}(f_{n})\|_{L^{2}\to L^{2}}\leq C \omega_{|t|}(nR)\|f_{n}\|_{\ell^{2}}.
\end{equation}
 For $R>0$ large enough and for any $n\in \mathbb{N}^{*}$, consider $f_{n}=\sum_{\gamma \in S^{\Gamma}_{n,R}}\mu_{n,R}(\gamma)D_{\gamma } \in \mathbb{C}[\Gamma]$ supported in $S^{\Gamma}_{n,R}$. Note that (\ref{growth}) of Theorem \ref{equi} implies the existence of some positive constant $C>0$ such that for any $n\in \mathbb{N}^*$ $$ \|f_{n}\|_{\ell^{2}}=\|\sum_{\gamma \in S^{\Gamma}_{n,R}}\mu_{n,R}(\gamma)D_{\gamma }\|_{\ell^{2}}\leq C/|S^{\Gamma}_{n,R}|^{\frac{1}{2}}.$$ 
 
 Lemma \ref{multiplicity} (\ref{item1}) implies the existence of $C>0$ such that  $|S^{\Gamma}_{n,R}|\leq C e^{-Q nR}$.
 From the lower bound of (\ref{HCHestimates}) together with above growth estimate we deduce the following ``spectral gap": there exists a constant $C>0$ such that for any $t\in \mathbb{R}$, we have for all $\gamma \in \Gamma,$ for all non negative integers $n$
\begin{equation}\label{radialestimates}
\|\pi_{t}(f_{n})\|_{L^{2}\to L^{2}}\leq C \widetilde{\phi_{t}}(nR),
\end{equation}

where $\widetilde{\phi_{t}}(nR)$ satisfies $C^{-1}\widetilde{\phi_{t}}(nR)\leq \phi_{t}(\gamma) \leq C\widetilde{\phi_{t}}(nR)$ for all $\gamma \in S^{\Gamma}_{n,R}$ where $C$ is a constant independent on $n$.\\ 
 The aim of this subsection is to prove a $L^{p}$-version of the above inequality (\ref{radialestimates}).\\


 

Although $\pi_{t}$ is an isometric action on $L^{p}$ with $1/p=1/2+t$, it defines also a representation $\pi_{t}:\Gamma \rightarrow  \mathbb{GL}(L^{r})$ where $ \mathbb{GL}(L^{r})$ stands for the group of bounded invertible linear operators acting on $L^{r}$. More precisely

\begin{prop}
For any $-1/2\leq t\leq 1/2$, for all $\gamma \in \Gamma$ the operator $\pi_{t}(\gamma)$ is bounded invertible operator on $L^{r}$ for all $1\leq r\leq \infty$ and moreover $\gamma \in \pi_{t}(\gamma) \in  \mathbb{GL}(L^{r})$ is a group morphism.
\end{prop}

\begin{proof}
Pick $-1/2\leq t\leq 1/2$.  Assume $1\leq r < \infty$. We have for all $\gamma \in \Gamma$ and for all $v\in L^{r}$
\begin{align*}
\|\pi_{t}(\gamma)v\|^{r}_{r}&=\int_{\partial X}e^{r(1/2+t)Q\beta_{\xi}(o,\gamma o)}|v(\gamma^{-1}\xi)|^{r}d\nu_{o}(\xi)\\
&=\int_{\partial X}e^{Q\beta_{\xi}(o,\gamma o) +(r(1/2+t)-1) Q\beta_{\xi}(o,\gamma o)  }|v(\gamma^{-1}\xi)|^{r}d\nu_{o}(\xi)\\
&\leq e^{|r(1/2+t)-1)| Q |\gamma|  } \int_{\partial X} e^{Q\beta_{\xi}(o,\gamma o) } |v(\gamma^{-1}\xi)|^{r}d\nu_{o}(\xi)\\
&= e^{|r(1/2+t)-1)| Q |\gamma|  } \|v\|^{r}_{r},
\end{align*}
where the inequality follows from the fact $|\beta_{\xi}(x,y)|\leq d(x,y)$ for all $\xi\in \partial X$ and for all $x,y\in X$. \\
 For the case $r=\infty$, we have  for all $\gamma \in \Gamma$ and for all $v\in L^{\infty}$ that $$\|\pi_{t}(\gamma)v\|_{\infty}\leq e^{|1/2+t| Q |\gamma|  } \|v\|_{\infty}.$$Hence, for all $\gamma \in \Gamma$ the operator $\pi_{t}(\gamma)$ is bounded invertible operator on $L^{r}$ for all $1\leq r\leq \infty$. The cocycle property of the Radon-Nikodym derivative implies that $\pi_{t}$ is a morphism and thus $\pi^{-1}_{t}(\gamma)=\pi_{t}(\gamma^{-1})$.  \end{proof}
In order to prove a $L^{p}$-version of  Inequality (\ref{radialestimates}) we need the following crucial lemma.\\
\begin{lemma}\label{BMtrick}
 Let $R>0$ large enough. For any $n\in \mathbb{N}$,  set $f_{n}=\sum_{\gamma \in S^{\Gamma}_{n,R}}\mu_{n,R}(\gamma)D_{\gamma } \in \mathbb{C}[\Gamma]$. Consider 
$\pi_{t}(f_{n}) $ as an operator from $L^{\infty} \to L^{\infty}$ with $t\in [-1/2,1/2]$.
The exists $C_{\infty}>0$ such that for all $t\in [-1/2,1/2]$ and for all $n\in \mathbb{N}$ $$ \|\pi_{t}(f_{n})\|_{L^{\infty} \to L^{\infty}}\leq C_{\infty} \widetilde{\phi_{t}}(nR).$$
\end{lemma}

\begin{proof}
Let $t\in \mathbb{R}$ such that $-1/2\leq t\leq1/2$ and $R>0$. Consider the sequence of functions $(G_{n})_{n}$ defined for each $n$ as $$G_{n}:\xi \in \partial X \mapsto \sum_{\gamma \in S^{\Gamma}_{n,R}}\mu_{n,R}(\gamma)e^{(\frac{1}{2}+t)Q\beta_{\xi}(o,\gamma o)}.$$ Consider also $$\check{G}_{n}:\xi \in \partial X \mapsto \sum_{\gamma \in S^{\Gamma}_{n,R}}\mu_{n,R}(\gamma^{-1})e^{(\frac{1}{2}+t)Q\beta_{\xi}(o,\gamma o)}. $$

Let $F\in L^{\infty}(\partial X,\nu_{o})$. We have for every $\xi\in \partial X$ and for all $n\in \mathbb{N}:$ 
\begin{align*}
|\pi_{t}(f_{n})F(\xi)|&=\big|\sum_{\gamma \in S^{\Gamma}_{n,R}}\mu_{n,R}(\gamma)e^{(\frac{1}{2}+t)Q\beta_{\xi}(o,\gamma o)}F(\gamma^{-1}\xi)\big|\\
&\leq \sum_{\gamma \in S^{\Gamma}_{n,R}}\mu_{n,R}(\gamma)e^{(\frac{1}{2}+t)Q\beta_{\xi}(o,\gamma o)}|F(\gamma^{-1}\xi)| \\
&\leq \sum_{\gamma \in S^{\Gamma}_{n,R}}\mu_{n,R}(\gamma)e^{(\frac{1}{2}+t)Q\beta_{\xi}(o,\gamma o)}\|F\|_{\infty},
\end{align*} 
and thus:
$$\|\pi_{t}(f_{n})F\|_{\infty}\leq  \|G_{n}\|_{\infty}\|F\|_{\infty}.$$
In other words, we shall prove $$\sup_{n}  \frac{ \|G_{n}\|_{\infty}}{\widetilde{\phi_{t}}(nR)}<\infty.$$ Pick $\xi\in \partial X$. Lemma \ref{crucial} implies that one can choose $R>0$ large enough such that there exist $C>0$ and $g_{\xi} \in S_{n,R}$ satisfying for all $\gamma \in S^{\Gamma}_{N,R}$
$$\beta_{\xi}(o,\gamma o)\leq  \beta_{\hat{\gamma}_{o}}(o,g_{\xi} o)+C.$$ Furthermore, the right inclusion of Lemma \ref{ombre} implies that there exists $C'>0$ such that for all $\eta \in O_{r}(o,\gamma o)$
$$  \beta_{\hat{\gamma}_{o}}(o,g_{\xi} o)\leq \beta_{\eta}(o,g_{\xi} o)+C'.$$
It follows a ``quasi mean-value property" that reads as follows 
$$ e^{(\frac{1}{2}+t)Q\beta_{\hat{\gamma}_{o}}(o,g_{\xi}o)} \leq \frac{C_{Q,t}}{|\nu_{o}(O_{r}(o,\gamma o))|} \int_{O_{r}(o,\gamma o)}e^{(\frac{1}{2}+t)Q\beta_{\eta }(o,g_{\xi}o)}d\nu_{o},$$
where $C_{Q,t}=e^{(\frac{1}{2}+t)Q+C'}.$
  Therefore, using an absorbing constant $C$ independent on $t\in [-1/2,1/2]$ we obtain

\begin{align*}
G_{n}(\xi)&\leq C \sum_{\gamma \in S^{\Gamma}_{n,R}}\mu_{n,R}(\gamma)e^{(\frac{1}{2}+t)Q\beta_{\hat{\gamma}_{o}}(o,g_{\xi}o)} \\
&\leq C \sum_{\gamma \in S^{\Gamma}_{n,R}} \frac{\mu_{n,R}(\gamma)}{\nu_{o}(O_{r}(o,\gamma o))}\int_{O_{r}(o,\gamma o)}e^{(\frac{1}{2}+t)Q\beta_{\eta}(o,g_{\xi}o)}d\nu_{o}(\eta) \\
&\leq  C \sum_{\gamma \in S^{\Gamma}_{n,R}} \int_{O_{r}(o,\gamma o)}e^{(\frac{1}{2}+t)Q\beta_{\eta }(o,g_{\xi}o)}d\nu_{o}(\eta) \\
&\leq  C \phi_{t}(g_{\xi}),
\end{align*}
where the first inequality follows Lemma \ref{crucial}, the second inequality follows from Lemma \ref{shadow} combined with the growth of $|S^{\Gamma}_{n,R}|$ and the last inequality follows from the finite multiplicity of the covering $\cup_{ \gamma \in S^{\Gamma}_{n,R} }O_{r}(o,\gamma o)$ proved in Lemma \ref{multiplicity}. \\
The estimates of the spherical functions (\ref{HCHestimates}) applied to $g_{\xi}\in S^{\Gamma}_{n,R}$ together with the above inequality imply for almost every $\xi\in 
\partial X$ and for all $n\in \mathbb{N}:$
$$G_{n}(\xi)\leq C \widetilde{\phi_{t}}(nR).$$ Thus
$$\sup_{n}\frac{\|G_{n}\|_{\infty}}{\widetilde{\phi_{t}}(nR)}<+\infty,$$ as required.

The above method applied to $\check{G}_{n}$ implies $$\sup_{n}\frac{\|\check{G}_{n}\|_{\infty}}{\widetilde{\phi_{t}}(nR)}<\infty.$$

\end{proof}
Eventually we obtain the $L^{p}$-version of Radial property RD for $L^p$-boundary representations.
\begin{theorem}\label{Lpspectral}

Let $R>0$ be large enough and let $r\in [1,\infty]$ such that $0\leq 1/r\leq 1$. There exists $C>0$ such that for any $-\frac{1}{2}\leq t \leq\frac{1}{2}$,  such that for all $n\in \mathbb{N}$ with $f_{n}=\sum_{\gamma \in S^{\Gamma}_{n,R}}\mu_{n,R}(\gamma)D_{\gamma } \in \mathbb{C}[\Gamma]$ supported on $S^{\Gamma}_{n,R}$ we have: 
 \begin{align}
  \|\pi_{t} (f_{n})\|_{L^{r}\to L^{r}}\leq C \widetilde{\phi_{t}}(nR).
  \end{align}
   And thus we have 
    \begin{align}
  \sup_{\|v\|_{p},\|w\|_{q}\leq 1} |\langle \pi_{t} (f_{n})v,w\rangle |\leq C\widetilde{\phi_{t}}(nR).
  \end{align}

 \end{theorem}
\begin{remark}
It is worth noting that the above theorem can be viewed as a $L^{p}$-version of Radial Property RD for $L^p$-boundary representations of hyperbolic groups.
\end{remark}
\begin{proof}

The proof is based on Riesz-Thorin Theorem. We shall prove that for any $t\in [-1/2,1/2 ] $ and for each $n\in \mathbb{N}$ the operator  $\pi_{t}(f_{n})$ viewed as  an operator from $L^{1} \to L^{1}$ and as an operator from $L^{\infty} \to L^{\infty}$ is uniformly bounded with respect to $n$.\\
The second point follows from Lemma \ref{BMtrick}.  We prove now the first point. First, observe that $\pi_{t}(f_{n})$ preserves the cone of positive functions since $f_{n}$ is positive and $\pi_{t}$ itself preserves the cone of positive functions. By decomposing a function $v$ into real an imaginary parts and positive and negative functions it is enough to find a  bound for positive functions. Assume $v \geq 0$.
\begin{align*}
\|\pi_{t}(f_{n})v\|_{1}&=\langle \pi_{t}(f_{n})v,\textbf{1}_{\partial X}\rangle \\
&=\langle v,  \pi_{-t}(\check{f_{n}}) \textbf{1}_{\partial X}\rangle \\
&\leq \|\check{G}_{n}\|_{\infty} \langle v,  \textbf{1}_{\partial X}\rangle \\
&\leq C_{\infty}\widetilde{\phi_{t}}(nR) \|v\|_{1},
\end{align*}
where the first inequality follows from the proof of Lemma \ref{BMtrick}.
Therefore Riesz-Thorin theorem implies that 
 $\pi_{t}(f_{n}) $ defines a bounded operator from $L^{r}$ to $L^{r}$ for any $r$ such that $1/r\in [0,1]$.
  \end{proof}

\section{Proofs}\label{section4}

The proof of Theorem \ref{BML2} relies on three steps. 
\begin{proof}
Let $0<t<\frac{1}{2}$. \\
\textbf{Step 1:} Uniform boundedness. Consider for $R>0$ and for all non-negative integer $n$ the function $f_{n}$ supported on $S^{\Gamma}_{n,R}$ defined as: $$f_{n}=\sum_{\gamma \in S^{\Gamma}_{n,R}}\mu_{n,R}(\gamma)\frac{ D_{\gamma }}{\phi_{t}(\gamma)}.$$ 

Note that this a weighted version of the function $\sum_{\gamma \in S^{\Gamma}_{n,R}}\mu_{n,R}(\gamma) D_{\gamma }$ by the spherical function $\phi_{t}(\gamma).$ The $L^{p}$-spectral inequality of Theorem \ref{Lpspectral} together with the fact that there exists $C>0$ such that for all $\gamma \in S^{\Gamma}_{n,R}$, $C^{-1}\widetilde{\phi_{t}}(nR)\leq \phi_{t}(\gamma)  $ imply
 
 \begin{align}
\sup_{n}\| \pi_{t}(f_{n})\|_{L^{p}\to L^{p}} <+\infty.
\end{align}

Set $K:=\sup_{n}\| \pi_{t}(f_{n})\|_{L^{p}\to L^{p}} $. Given $f,g\in C(\overline{ X})$ we have for all $n\in \mathbb{N}$

 \begin{align*}
| \sum_{\gamma \in S^{\Gamma}_{n,R}}\mu_{n,R}(\gamma)f(\gamma  o)g(\gamma^{-1} o)\frac{\langle\pi_{t}(\gamma)v,w\rangle }{\phi_{t}(\gamma)}|&\leq C \|f\|_{\infty}\|g\|_{\infty}\| \pi_{t}(f_{n})\|_{L^{p}\to L^{p}} \|v\|_{p}   \|w\|_{q} \\
 &\leq K \|f\|_{\infty}\|g\|_{\infty} \|v\|_{p}   \|w\|_{q}.
 \end{align*}

By Banach-Alaoglu-Bourbaki Theorem and since on reflexive spaces the weak topology and the weak*-topology coincide,     the limit $$\lim_{n\to +\infty} \sum_{\gamma \in S^{\Gamma}_{n,R}}\mu_{n,R}(\gamma)f(\gamma  o)g(\gamma^{-1} o)\frac{\langle\pi_{t}(\gamma)v,w\rangle }{\phi_{t}(\gamma)}$$ exists for all $v\in L^{p}$, for all $w\in L^{q}$ and for all $f,g\in C(\overline{X})$, up to extraction.\\

\textbf{Step 2:} Computation of the limit.\\

 We already know by \cite{BPi} that we have the desired result for $v,w$ in a dense subspace of $L^{p},L^{q}$ spaces (e.g. that for all $v,w \in Lip(\partial X)$ and $f,g\in C(\overline{X})$)

\begin{align*}
 \lim_{n\to +\infty}\sum_{\gamma \in S^{\Gamma}_{n,R}}\mu_{n,R}(\gamma)f(\gamma  o)g(\gamma^{-1} o)\frac{\langle\pi_{t}(\gamma)v,w\rangle }{\phi_{t}(\gamma)}=\langle g_{|_{\partial X}} \mathcal{R}_{t}(v), \textbf{1}_{\partial X}\rangle  \langle  f_{|_{\partial X}} ,w \rangle.
 \end{align*}

\textbf{Step 3:} Conclusion.\\
Assume that $0<t<1/2$. Therefore $\mathcal{I}_{t}$ and thus $\mathcal{R}_{t}$ are continuous on $L^{p}$ with $1/p=1/2+t$.
 The limit above together with the uniform bound of \textbf{Step 1} imply eventually that  for $f,g\in C(\overline{ X})$ and for 
 $(v,w)\in L^{p}\times L^{q}$: 
\begin{align*}
 \lim_{n\to +\infty}\sum_{\gamma \in S^{\Gamma}_{n,R}}\mu_{n,R}(\gamma)f(\gamma  o)g(\gamma^{-1} o)\frac{\langle\pi_{t}(\gamma)v,w\rangle }{\phi_{t}(\gamma)}=\langle g_{|_{\partial X}} \mathcal{R}_{t}(v), \textbf{1}_{\partial X}\rangle \overline{ \langle  w,f_{|_{\partial X}}  \rangle}.
 \end{align*}

\end{proof}

\subsection{Proof of irreducibility}

To prove irreducibility of representations our main tool is Theorem \ref{BML2}. 
  \begin{lemma}\label{cyclic}
Let $0<t<\frac{1}{2}$. The vector $\textbf{1}_{\partial X}$ is cyclic for $(\pi_{t},L^{p})$ and the vector $\sigma_{t}$ is cyclic for  $(\pi_{-t},\overline{Im(\mathcal{I}_{t})} ^{\|\cdot\|_{q}})$.
  \end{lemma}
  
  \begin{proof}
Theorem  \ref{BML2} implies for all $w\in L^{q}$ the following convergence and for all continuous functions $f\in C(\overline{X})$:
  
$$\lim_{n\to +\infty}\sum_{\gamma \in S^{\Gamma}_{n,R}}\mu_{n,R}(\gamma)f(\gamma  o)\frac{\langle\pi_{t}(\gamma)\textbf{1}_{\partial X},w\rangle}{\phi_{t}(\gamma)}=  \langle  f_{|_{\partial X}},w \rangle. $$
Now, given a function $f\in C(\partial X)$ consider $\tilde{f}\in C(\overline{X})$ such that $\tilde{f}_{|_{\partial X}}=f$. Set 
\begin{equation}\label{fn}
f_{n}:=\sum_{\gamma \in S^{\Gamma}_{n,R}}\mu_{n,R}(\gamma)\tilde{f}(\gamma  o)\frac{\pi_{t}(\gamma)\textbf{1}_{\partial X}}{\phi_{t}(\gamma)}\in \pi_{t}(\mathbb{C}[\Gamma])\textbf{1}_{\partial X}.
\end{equation}
Hence, $f_{n}\to f$ with respect to the weak topology of $L^{p}$. Therefore, since $C(\partial X)$ is dense in $L^{p}$ and since the closure of the weak topology and the $\|\cdot\|_{p}$ coincides, the vector $\textbf{1}_{\partial X}$ is cyclic for $L^{p}$.\\ 
We prove now that $\sigma_{t}$ is cyclic for $\pi_{-t}$. Recall that for all $v\in L^{p}$ we have $\mathcal{I}_{t}(v)\in L^{q}$ by Proposition \ref{cont}. Hence for all $v\in L^{p}$, with the same notation of (\ref{fn}) we have

$$\langle f_{n},\mathcal{I}_{t}(v) \rangle \to \langle  f_{|_{\partial X}},\mathcal{I}_{t}(v)\rangle.$$ Since $\mathcal{I}_{t}$ is self adjoint, for all $v\in L^{p}$ 
\begin{align*}
\langle \mathcal{I}_{t}(f_{n}),v \rangle&=\langle \sum_{\gamma \in S^{\Gamma}_{n,R}}\mu_{n,R}(\gamma)\tilde{f}(\gamma  o)\frac{\pi_{-t}(\gamma)\mathcal{I}_{t}(\textbf{1}_{\partial X})}{\phi_{t}(\gamma)},v\rangle\\
&\to \langle \mathcal{I}_{t}(f_{|_{\partial X}}),v \rangle,\\
\end{align*}
where $ \mathcal{I}_{t}(f_{n}) \in \pi_{-t}(\mathbb{C}[\Gamma])\sigma_{t}.$ 

Hence $$\overline{ \pi_{-t}(\mathbb{C}[\Gamma])\sigma_{t}} ^{\|\cdot\|_{q}}=\overline{\mathcal{I}_{t}(C(\partial X))}^{\|\cdot\|_{q}}.$$

Eventually, using the density of $C(\partial X)$ in $L^{p}$ and the continuity of $\mathcal{I}_{t}$ we deduce $$\overline{ \pi_{-t}(\mathbb{C}[\Gamma])\sigma_{t}} ^{\|\cdot\|_{q}}= \overline{\mathcal{I}_{t}(L^{p})}^{\|\cdot\|_{q}}.$$


\end{proof}


\begin{proof}[Proof of Theorem \ref{mainT}.]

First of all, since $\mathcal{I}_{t}$ is continuous from $L^{p}$ to $L^{q}$, the subspace $\ker \mathcal{I}_{t}$ is a closed invariant subspace of $(\pi_{t},L^{p})$. Thus, if $\mathcal{I}_{t}$ is non injective, then $(\pi_{t},L^{p})$ is not irreducible.\\
We shall prove now that $\mathcal{I}_{t}$ is injective then $(\pi_{t},L^{p})$ is irreducible for $1/p=1/2+t.$ Since $\mathcal{I}_{t}$ is a continuous operator from $L^{p}\to L^{q}$, a standard result in Banach spaces theory asserts that the dual space of $L^{p}/ \ker \mathcal{I}_{t}$ is the space $\overline{ Im (\mathcal{I}_{t})}^{\|\cdot\|_{q}}$, where the  the weak closure for $L^{q}$ is the same as the $\|\cdot\|_{q}$-closure (see \cite[Chapter 3]{Me}). Hence, the dual representation of $(\pi_{t}, L^{p}/ \ker \mathcal{I}_{t} )$ is $(\pi_{-t}, \overline{ Im (\mathcal{I}_{t})}^{\|\cdot\|_{q}})$.  We shall prove that $(\pi_{-t},L^{q})$ is irreducible to obtain irreducibility of $(\pi_{t}, L^{p}/ \ker \mathcal{I}_{t} )$.\\ Recall that 
 $ \pi_{-t}(f) \sigma_{t}=\mathcal{I}_{t}(\pi_{t}(f)\textbf{1}_{\partial X} )$, for all $f\in \mathbb{C}[\Gamma]$ and the function $\sigma_{t}$ is cyclic in $\overline{ Im (\mathcal{I}_{t})}^{\|\cdot\|_{q}}$ for $\pi_{-t}$ by Lemma \ref{cyclic}.\\
 Now, let $0\neq K\subset \overline{ Im (\mathcal{I}_{t})}^{\|\cdot\|_{q}} \subset L^{q}(\partial X,\nu_{o})$ a $L^{q}$-closed subspace invariant by $\pi_{-t}$. Let $R>0$ large enough. For any $w\in K \subset L^{q}$ define for all $n\in \mathbb{N}$, the vector:
 \begin{equation}\label{zzz}
  w_{n}:=\sum_{\gamma \in S^{\Gamma}_{n,R}}\mu_{n,R}(\gamma) \tilde{ \sigma_{t}}(\gamma^{-1} o)  \frac{\pi_{-t}(\gamma^{-1}) w }{\phi_{t}(\gamma)} \in K,
  \end{equation}
  
  where $\tilde{\sigma_{t}}$ has been defined in (\ref{sigmat}).
  Theorem \ref{BML2} implies for all $v\in L^{p}(\partial X,\nu_{o})$ that as $n\to +\infty$:

$$ \sum_{\gamma \in S^{\Gamma}_{n,R}}\mu_{n,R}(\gamma) \tilde{\sigma_{t}}(\gamma^{-1} o) \frac{\langle  v, \pi_{-t}(\gamma^{-1}) w\rangle }{\phi_{t}(\gamma)}\to \langle \mathcal{I}_{t}(v),\textbf{1}_{\partial X}\rangle \langle  \textbf{1}_{\partial X} ,w \rangle=\langle v,\mathcal{I}_{t}(\textbf{1}_{\partial X})\rangle \langle\textbf{1}_{\partial X}, w \rangle,$$

and the above convergence reads as follows with respect to the weak topology on $L^{q}:$
$$w_{n}\to   \langle \textbf{1}_{\partial X} , w \rangle \sigma_{t}.$$
Since $K$ is closed we have that $ \langle \textbf{1}_{\partial X} , w \rangle \sigma_{t} \in K.$ So, since $\sigma_{t}$ is cyclic, it is sufficient to show that there exists $0\neq w \in K$ such that $$\langle \textbf{1}_{\partial X} , w \rangle\neq 0.$$ 

Assume this is not the case: for all $w\in K$ we have $\langle \textbf{1}_{\partial X},w\rangle =0$. We would have that   for all $\gamma \in \Gamma$ that $\langle  \textbf{1}_{\partial X},\pi_{-t}(\gamma^{-1})w \rangle=0= \langle \pi_{t}(\gamma) \textbf{1}_{\partial X},w \rangle.$ And therefore, since $\textbf{1}_{\partial X}$ is cyclic for $(\pi_{t},L^{p}(\partial X,\nu_{o}))$ it implies that
$$  K\ \subset \{w\in L^{q}| \langle v,w\rangle=0,\forall v\in L^{p} \}$$

Since the pairing is non-degenerate then $ K=\{0 \}.$

Hence, if $K\neq \{0 \}$ then $K$ contains $\sigma_{t}$ and thus it has to be $\overline{ Im (\mathcal{I}_{t})}^{\|\cdot\|_{q}}$ and the proof is done.

\end{proof}
\section{Application to rank one semisimple Lie groups}\label{section5}
Let $G$ be a connected semisimple Lie group with finite center and let $\frak g$ be its Lie algebra.
Let $K$ be a maximal compact subgroup of $G$ and let $\frak k$ be its Lie algebra. Let $\frak p$ be the orthogonal complement of $\frak k$ in $\frak g$
relative to the Killing form $B$. Among the  abelian sub-algebras of $\frak g$ contained in the subspace $\frak p$, let $\frak a$ be a maximal one. We assume $\dim \frak a=1$, i.e. the real rank of $G$ equals to $1$ (in particular $G$ is not compact). Let $\Sigma\subset \frak{a}^{*}$ be the root system associated to $(\frak g,\frak a)$. Let
\[
	\frak g_{\alpha}=\{X\in\frak g: \mathrm{ad}(H)X=\alpha(H)X\ \ \forall H\in\frak a  \}
\]
be the root space of $\alpha\in\Sigma$. Recall that $\Sigma=\{-\alpha,\alpha\}$ or $\Sigma=\{-2\alpha,-\alpha,\alpha,2\alpha\}$ where $\alpha$ is a positive root ($\alpha\in\Sigma$ is positive if and only if $\alpha(H)>0$ for all $H\in \frak a^+$). If $m_{1}=\dim \frak g_{\alpha}$ and $m_{2}=\frak g_{2\alpha}$ denoted by $\rho=\frac{1}{2}(m_{1}+2m_{2})\alpha$ be the half sum of positive roots. Let $H_{0}$
 be the unique vecteur in $\frak{a}$ such that $\alpha(H_{0})=1$. Hence, $\frak{a}=\{ tH_{0},t\in\mathbb{R}\}$ and $\frak{a}_{+}=\{ H\in \frak{a}|\alpha(H)>0\}$ is identified with the open set of strictly positive real numbers.
Let $\frak n$ the nilpotent Lie algebra defined as the direct sum of root spaces of positive roots:
\[
	\frak n=\bigoplus_{\alpha\in\Sigma^+}\frak g_{\alpha}.
\]
Let $A=\exp(\frak a)$, $A^+=\exp(\frak a^+)$ and $N=\exp(\frak n)$.
Let $G=KAN$ be the Iwasawa decomposition and $K\overline{A^{+}}K$ the Cartan decomposition defined by $\frak a^+$ where $\overline{A^{+}}$ denoted the closure $A^{+}$. Let $Z(A)$ be the centralizer of $A$ in $G$ and  
 $M=Z(A)\cap K$. The  group $M$ normalizes $N$. Let $P=MAN$ be the minimal parabolic subgroup of $G$ associated to $\frak a^+$.
Let $\nu$ be the unique Borel regular $K$-invariant probability measure on the Furstenberg-Poisson boundary $G/P$  that is quasi-invariant under the action $G$ (we refer to \cite[Appendix B]{BDV} for a general discussion). 
Let 
\[
\rho_{t}:G\to \mathcal{U}(L^p(G/P,\nu))	
\]
be the associated $L^{p}$-boundary representation of $G$ and defined the corresponding spherical function
\[
	\phi_{t}(g)=\langle \rho_{t}(g)\textbf{ 1}_{G/P}, \textbf{1}_{G/P}\rangle.\]

The corresponding globally symmetric space of non compact type of $G$ is  $G/K$  endowed with a $G$-invariant Riemannian metric denoted by $d$ induced by the  Killing form on $\frak{g}\ / \frak{k}$ identified with the tangent space of $G/K$ at the point $o=eK$.  A flat of dimension $k$ is defined as the image of a map $\mathbb{R}^{k}\rightarrow G/K$ locally isometric.  The rank of $G$ is the largest dimension of a flat subspace of $G/K$.  The rank one globally symmetric spaces of non compact type are classified as follows: there are  the real hyperbolic spaces $H^{n}(\mathbb{R})$, the complex hyperbolic spaces $H^{n}(\mathbb{C})$, the quaternionic hyperbolic spaces $H^{n}(\mathbb{H})$ for $n\geq 2$ and the exceptional hyperbolic space that is the $2$-dimensional octonionic hyperbolic space $H^{2}(\mathbb{O})$. \\
If $(X,d)$ is one the above hyperbolic space, it is a CAT(-1) space and in particular it is a proper geodesic $\delta$-hyperbolic space and fits of course in the class of spaces \ref{class} for $\epsilon=1$. One can therefore consider its Gromov boundary $\partial X$ or equivalently the geometric boundary of $X$.  The group $G$ acts by isometries on $(X,d)$ and its discrete subgroup acts properly discontinuously on $X$. Assume that $\Gamma$ is lattice (uniform or non-uniform) and perform the Patterson-Sullivan construction associated to $(\Gamma,d)$ with the base point $o=eK \in X$ to obtain a measure supported on $\partial X$ denoted by $\nu_{o}$. The Hausdorff dimension of $\nu_{o}$ is the critical exponent  $Q_{\Gamma}$ of $\Gamma$ that coincides with the volume growth of the corresponding hyperbolic spaces $Q_{G}=m_{1}+2m_{2}$.\\
 A geodesic ray starting at the origin can be represented  using Cartan decomposition as $c(t)=ke^{tH_{0}}\cdot o=kMe^{tH_{0}}\cdot o$ where $t\in \mathbb{R}_{+}$ and $k\in K$. Then the Furstenberg-Poisson boundary $G/P$ can be identified with the geometric boundary $\partial X$ in the case of rank one symmetric space. Indeed, one can identify $G/P=K/M$ thanks to the Iwasawa decomposition $KAN$. It turns out that the Patterson-Sullivan measure $\nu_{o}$ associated with a lattice $\Gamma$ supported on $\partial X$ of dimension $Q_{\Gamma}$ coincides with the unique $K$-invariant measure $\nu$ on $G/P$. \\
 Thus, the $L^{p}$-boundary representation of $\Gamma$ is nothing but the restriction of $\rho_{t}$ to $\Gamma$
 \begin{align}
 \pi_{t}:\gamma \in \Gamma\rightarrow \rho_{t_{|_{\Gamma}}}   \in \mbox{Iso}(L^{p}(\partial X,\nu_{o}))
 \end{align} 
 with $p$ such that $1/p=1/2+t$ with $0<t<1/2$. Since $\Gamma$ might be a non-uniform lattice, the results obtained above dealing with hyperbolic groups do not apply to $\Gamma$. Nevertheless we have the exact analog of Lemma \ref{BMtrick}.\\
In the following, for $R>0$ and for any $n\in \mathbb{N}$, the spheres $S^{\Gamma}_{n,R}$  are defined with respect to the length function $|\gamma|:=d(o,\gamma o)$ corresponding to the Riemannian metric $d$ on the symmetric space $G/K$. Moreover, one can take for $\mu_{n,R}$ , for $R>0$ and for any $n\in \mathbb{N}$ the standard average
$$\frac{1}{|S^{\Gamma}_{n,R}|} \sum_{\gamma \in S^{\Gamma}_{n,R}} D_{\gamma }.$$
\begin{lemma}\label{BMtrick2}
Let $\Gamma$ be a lattice in $G$.
 Let $R>0$.  For all $t\in [-1/2,1/2]$ and  $n\in \mathbb{N}$ set $f_{n}=\frac{1}{|S^{\Gamma}_{n,R}|}\sum_{\gamma \in S^{\Gamma}_{n,R}}D_{\gamma } \in \mathbb{C}[\Gamma]$. Consider 
$\pi_{t}(f_{n}) $  as an operator from $L^{\infty} \to L^{\infty}.$
The exists $C_{\infty}>0$ depending on $R$  such that for all $t\in [-1/2,1/2]$ $$ \|\pi_{t}(f_{n})\|_{L^{\infty} \to L^{\infty}}\leq C_{\infty} \widetilde{\phi_{t}}(nR).$$
\end{lemma} 

\begin{proof}
Indeed the proof follows the same ideas of \cite[Proposition 3.2]{BoyHCH} and \cite[Section 2.5]{BLP}
\end{proof}
Therefore, following exactly the same method of the proof of Theorem \ref{BML2}, we obtain

\begin{theorem}\label{lattspectral}
Let $\Gamma$ be a lattice in a rank one connected semisimple Lie group with finite center $G$. 
Let $R>0$ be large enough and let $r\in [1,\infty]$ such that $0\leq 1/r\leq 1$. There exists $C>0$ such that for any $-\frac{1}{2}\leq t \leq\frac{1}{2}$,  such that for all $n\in \mathbb{N}$ with $f_{n}=\frac{1}{|S^{\Gamma}_{n,R}|}\sum_{\gamma \in S^{\Gamma}_{n,R}}D_{\gamma } \in \mathbb{C}[\Gamma]$ supported on $S^{\Gamma}_{n,R}$ we have: 
 \begin{align}
  \|\pi_{t} (f_{n})\|_{L^{r}\to L^{r}}\leq C \widetilde{\phi_{t}}(nR).
  \end{align}
   And thus we have 
    \begin{align}
  \sup_{\|v\|_{p},\|w\|_{q}\leq 1} |\langle \pi_{t} (f_{n})v,w\rangle |\leq C\widetilde{\phi_{t}}(nR).
  \end{align}
\end{theorem}

The equidistribution theorem needed in (\ref{equid}) reads as follows in the context of lattices. We refer to \cite{Ro} for the next results in a more general setting. 

\begin{theorem}\label{equi}
Let $\Gamma$ be a lattice in a rank one connected semisimple Lie group with finite center $G$. For any  $R>0$, we have the following convergence: $$\frac{1}{|S^{\Gamma}_{n,R}|}\sum_{\gamma \in S^{\Gamma}_{n,R}} D_{\gamma o} \otimes D_{\gamma^{-1} o} \rightharpoonup  \nu_{o}\otimes \nu_{o},$$ 
as $n\to +\infty$, for the weak* convergence in $C(\overline{X}\times \overline{X})$.
\end{theorem}

Eventually we obtain
\begin{theorem}\label{BML3}
Let $\Gamma$ be a lattice in a rank one connected semisimple Lie group with finite center $G$. For  $R>0$ large enough, for all $0<t<1/2  $,  for all $f,g\in C(\overline{X})$, for all $v\in L^{p}( \partial X,\nu_{o})$ and $w\in L^{q}(\partial X,\nu_{o})$:
$$\frac{1}{|S^{\Gamma}_{n,R}|}\sum_{\gamma \in S^{\Gamma}_{n,R}} f(\gamma   ) g(\gamma^{-1} ) \frac{\langle \pi_{t}(\gamma)v,w\rangle }{\phi_{t}(\gamma)}\to \langle  g_{|_{\partial \Gamma}}\mathcal{R}_{t}(v),\textbf{1}_{\partial \Gamma}\rangle \langle  f_{|_{\partial \Gamma}},w  \rangle, $$
as $n\to +\infty$.
\end{theorem}

\subsection{Intertwining operators}
The restriction of $\rho_{t}$ to $K$ on $L^{2}(K/M,\nu_{o})$ is given by 
\begin{align*}
\rho_{t}(k)v(\xi)=v(k^{-1}\xi),
\end{align*}
since $\nu_{o}$ is $K$-invariant and note that the representation does not depend on $t$ anymore and provides a unitary representation of $K$. Therefore the intertwining relations \ref{intert} reads as follows:  for all $k\in K$ and for all $t>0$
$$ \mathcal{I}_{t}\rho_{t}(k)= \rho_{t}(k) \mathcal{I}_{t}.$$
Peter-Weyl Theorem implies that $$ L^{2}(K/M,\nu_{o})=\oplus_{n\geq 0} V_{n},$$ where $V_{n}$ are finite dimensional irreducible unitary representation of $K$. Therefore, Schur's lemma implies that there exists a sequence of scalars $(\lambda_{n})_{n\geq 0}$ such that $\mathcal{I}_{t}$ restricted to $V_{n}$ is a scalar operator as follows $ \mathcal{I}_{t |_{V_{n}}}=\lambda_{n} Id_{ |_{V_{n} } } $ with $\lambda_{n}\neq 0$ for all $n\geq 0$. We deduce that $\mathcal{I}_{t}$ is injective viewed as an operator acting on $L^{2}$ and therefore  it is injective as an operator from $L^{p}$ to $L^{q}$. Apply Theorem \ref{BML3} to lattices in rank one semisimple Lie groups and use exactly the same arguments of Theorem \ref{mainT} to complete the proof of Theorem \ref{latt}.


\begin{thebibliography}{99} 


\bibitem{Ank}
 J.-P. Anker, \emph{La forme exacte de l'estimation fondamentale de Harish-Chandra. (French) [The exact
form of Harish-Chandra's fundamental estimate]} C. R. Acad. Sci. Paris Sr. I Math. 305 (1987), no.
9, 371-374. 




\bibitem{BM}
U. Bader, R. Muchnik, \emph{Boundary unitary representations, irreducibility and rigidity.} J. Mod. Dyn.
5 (2011), no. 1, 49-69.

\bibitem{BM2}
U. Bader, J. Dymara, \emph{Boundary unitary representations right-angled hyperbolic buildings.} J. Mod. Dyn. 10 (2016), 413-437.







\bibitem{BHM}
S. Blach\`ere,  P. Ha\"issinsky, P. Mathieu, \emph{Harmonic measures versus quasiconformal measures for hyperbolic groups.} Ann. Sci. Ec. Norm. Sup\'er. 44, no. 4 (2011), 683 -- 721.




\bibitem{Bou}
K. Boucher, \emph{Analogs of complementary series for CAT(-1) groups},
	arXiv:2007.15369, 2020.



\bibitem{Boy}
A. Boyer, \emph{Some spherical functions on hyperbolic groups}, Journal of topology and analysis, 2019.
\bibitem{Boy2}
 
A. Boyer, \emph{Equidistribution, ergodicity, and irreducibility in CAT(-1) spaces.} Groups Geom. Dyn. 11 (2017), no. 3, 777-818.
\bibitem{BoyMa}
A. Boyer, D. Mayeda \emph{Equidistribution, ergodicity and irreducibility associated with Gibbs measures.} Comment. Math. Helv. 92 (2017), no. 2, 349-387.


\bibitem{BoyHCH}
A. Boyer \emph{Harish-Chandra‘s Schwartz algebras associated with discrete subgroups of semisimple Lie groups}, J. Lie Theory 27 (2017), no. 3, 831-844.


\bibitem{BPi}
A. Boyer and J-C. Picaud, \emph{Riesz operator and some spherical representations for hyperbolic groups}, 	arXiv:2201.00077, 2021.
 



\bibitem{BoyP}
 A. Boyer and A. Pinochet-Lobos, \emph{An ergodic theorem for the quasi-regular representation of the free
group.} Bull. Belg. Math. Soc. Simon Stevin 24 (2017), no. 2, 243-255.
\bibitem{BLP}
A. Boyer, G. Link, Ch. Pittet, \emph{Ergodic boundary representations.} Ergodic Theory Dynam. Systems 39 (2019), no. 8, 2017-2047.

\bibitem{BG}
A. Boyer and L. Garncarek, \emph{Asymptotic Schur orthogonality in hyperbolic groups with application to monotony},  Trans. Amer. Math. Soc. 371 (2019), 6815-6841.
 
 
 \bibitem{BH}	
 M.R. Bridson; A. Haefliger, \emph{Metric spaces of non-positive curvature.} Grundlehren der Mathematischen Wissenschaften [Fundamental Principles of Mathematical Sciences], 319. Springer-Verlag, Berlin, 1999. 

 
  \bibitem{BMo}
 M. Burger and S. Mozes, \emph{CAT(-1)-spaces, divergence groups and their commensurators,} J. Amer. Math. Soc. 9 (1996), no. 1, 57-93,
 
 \bibitem{Ca}
Pierre-Emmanuel Caprace, Mehrdad Kalantar, Nicolas Monod, \emph{ A type I conjecture and boundary representations of hyperbolic groups}, arXiv:2110.00190, 2021.


\bibitem{Co}
M. Coornaert, \emph{Mesures de Patterson-Sullivan sur le bord d'un espace hyperbolique au sens de Gromov (French, with French summary).} Pacific J. Math. 159 (1993), no. 2, 241-270.

\bibitem{Cow}
F. Astengo, M. Cowling, and B. Di Blasio, \emph{The Cayley transform and uniformly
bounded representations}, J. Funct. Anal., vol. 213, no. 2, pp. 241–269, 2004.

\bibitem{CM}
C. Connell, R. Muchnik, \emph{Harmonicity of quasiconformal measures and Poisson boundaries of hyperbolic spaces.} Geom. Funct. Anal. 17 (2007), no. 3,








\bibitem{FiPi}
A. Fig\`a-Talamanca, M.A. Picardello \emph{Harmonic analysis on free groups,} Lecture notes in pure and applied mathematics vol. 87, 1983.




\bibitem{Fink}
V. Finkelshtein, \emph{Diophantine Properties of Groups of Toral Automorphisms.}
Thesis (Ph.D.) University of Illinois at Chicago. 2017. 77 pp.

\bibitem{Ga}
  L. Garncarek, \emph{Boundary representations of hyperbolic groups},  arXiv:1404.0903, 2014.
 
 
 
  \bibitem{G}
E. Ghys, P. de la Harpe. \emph{Panorama. (French) Sur les groupes hyperboliques d'apr\`es Mikhael Gromov }(Bern, 1988), 1-25, Progr. Math., 83, Birkh\"auser Boston, Boston, MA, 1990.

  \bibitem{Gr}
  M. Gromov, \emph{Hyperbolic groups}, Essays in group theory, Springer, coll.``MSRI Publ.'' (no 8), 1987, p. 75-263.
  
	
New York, 2001.


\bibitem{GAG}
G-A. Gruetzner, \emph{Harmonic analysis on the boundary of hyperbolic groups}, preprint, https://hal.inria.fr/hal-03963498/,  2023.





\bibitem{KS}
M-G. Kuhn, T. Steger, \emph{More irreducible boundary representations of free groups.} Duke Math. J. 82
(1996), no. 2, 381-436.

\bibitem{KS2}
M-G. Kuhn, T, Steger, \emph{Monotony of certain free group representations.} J. Funct. Anal. 179 (2001),
no. 1, 1-17.
\bibitem{Ma}
G. Margulis, \emph{On some aspects of the theory of anosov systems.} Springer Monographs in Mathematics, Springer-Verlag, Berlin, With a survey by Richard Sharp: Periodic orbits of hyperbolic flows,
Translated from the Russian by Valentina Vladimirovna Szulikowska, 2004.


	\bibitem{Ro}
	T. Roblin, \emph{Ergodicit\'e et unique ergodicit\'e du feuilletage horosph\'erique, m\'elange du flot g\'eod\'esique et \'equidistributions diverses dans les groupes discrets en courbure n\'egative}, M\'emoires de la Soci\'et\'e Math\'ematique de France, Nouvelle S\'erie 95 (2003).
		\bibitem{Me}
		
		Robert E.  Megginson,  \emph{An introduction to Banach space theory.} Graduate Texts in Mathematics, 183. Springer-Verlag, New York, 1998. 
			
		
 \bibitem{NS}
     B. Nica and J.  \v{S}pakula, \emph{Strong hyperbolicity.} Groups, Geometry and Dynamics Volume 10, Issue 3, 2016, pp. 951-964. 
\bibitem{Pa}
S. J. Patterson, \emph{The limit set of a Fuchsian group.} Acta Math. 136 (1976), no. 3-4, 241-273,

		\bibitem{Su}
		D. Sullivan, \emph{The density at infinity of a discrete group of hyperbolic motions}, Int. Hautes Etudes Sci. Publ. Math. 50 (1979), 171-202.
		
				 
		 \bibitem{Tao}
T. Tao, \emph{Lecture Notes 2}, \href{https://www.math.ucla.edu/~tao/247a.1.06f/notes2.pdf}{Notes 247 A}.







   \end{thebibliography}
\end{document}